\documentclass[12pt,a4paper,onecolumn]{article}
\usepackage{textcomp}
\usepackage{tikz}
\usepackage{amssymb}
\usepackage{caption}
\usepackage{color}
\usepackage{verbatim}
\usepackage{fancybox}
\usepackage{titlesec} 
\renewcommand\thesection{\arabic{section}} 
\usepackage{mathtools}
\usepackage{amsmath}
\usepackage{mathrsfs}
\usetikzlibrary{positioning,fit,calc}
\usepackage[a4paper,left=1in,top=1in,right=1in,bottom=1in,nohead]{geometry} 
\usepackage[hidelinks]{hyperref}
\usepackage{array}
\usepackage{bm}
\usepackage[symbol]{footmisc} 
\def\correspondingauthor{\footnote{Corresponding author. Email: williewong088@gmail.com.}}

\tikzset{block/.style={draw,thick,text width=2cm,minimum height=1cm,align=center},
         line/.style={-latex}}
\newcolumntype{P}[1]{>{\centering\arraybackslash}m{#1}} 

\titleformat{\section}[block]{\large\scshape\bfseries}{\thesection.}{1em}{} 
\titleformat{\subsection}[block]{\bfseries}{\thesubsection.}{1em}{} 

\newtheorem{defn}{Definition}[section]
\newtheorem{thm}[defn]{Theorem}

\newtheorem{ppn}[defn]{Proposition}
\newtheorem{rmk}[defn]{Remark}

\newtheorem{cor}[defn]{Corollary}
\newtheorem{lem}[defn]{Lemma}

\begin{document}
\thispagestyle{empty}
\pagenumbering{arabic}
\begin{center}
    \textbf{\Large On optimal orientations of complete tripartite graphs}
\vspace{0.1 in} 
    \\{\large W.H.W. Wong\correspondingauthor{}, E.G. Tay}
\vspace{0.1 in} 
\\National Institute of Education, Nanyang Technological University, Singapore
\end{center}

\begin{abstract}
\indent Given a connected and bridgeless graph $G$, let $\mathscr{D}(G)$ be the family of strong orientations of $G$. The orientation number of $G$ is defined to be $\bar{d}(G):=min\{d(D)|D\in \mathscr{D}(G)\}$, where $d(D)$ is the diameter of the digraph $D$. In this paper, we focus on the orientation number of complete tripartite graphs. We prove a conjecture raised by Rajasekaran and Sampathkumar \cite{RG SR}. Specifically, for $q\ge p\ge 3$, if $\bar{d}(K(2,p,q))=2$, then $q\le{{p}\choose{\lfloor{p/2}\rfloor}}$. We also present some sufficient conditions on $p$ and $q$ for $\bar{d}(K(p,p,q))=2$.
\end{abstract}

\section{Introduction}
\indent\par Let $G$ be a simple connected graph with vertex set $V(G)$ and edge set $E(G)$. For any vertex $v\in V(G)$, its $\textit{eccentricity}$ $e(v)$ is defined as $e(v):=max\{d_G(v,x)|x\in V(G)\}$, where $d_G(v,x)$ is the length of a shortest $v-x$ path. The $\textit{diameter}$ of $G$, denoted by $d(G)$, is defined as $d(G):=max\{e(v)|v\in V(G)\}$. For a digraph $D$, the above notations are defined similarly. The \textit{outset} and \textit{inset} of a vertex $v\in V(D)$ are defined to be $O_D(v):=\{x\in V(D)|\text{ } v\rightarrow x\}$ and $I_D(v):=\{y\in V(D)|\text{ } y\rightarrow v\}$ respectively. If there is no ambiguity, we shall omit the subscript for the above notations.
\indent\par An $\textit{orientation}$ of a graph $G$ is a digraph obtained from $G$ by assigning to each edge $e\in E(G)$ a direction. An orientation $D$ of $G$ is \textit{strong} if every two vertices in $V(D)$ are mutually reachable. An edge $e\in E(G)$ is a \textit{bridge} if $G-e$ is disconnected. Robbins' well-known one-way street theorem [10] states that a connected graph $G$ has a strong orientation if and only if no edge of $G$ is a bridge.
\indent\par Given a connected and bridgeless graph $G$, let $\mathscr{D}(G)$ be the family of strong orientations of $G$. The \textit{orientation number} of $G$ is defined to be $\bar{d}(G):=min\{d(D)|D\in \mathscr{D}(G)\}$, where $d(D)$ is the diameter of the digraph $D$. Trivially, $d(D)\ge d(G)$ for any $D\in \mathscr{D}(G)$. An orientation $D\in \mathscr{D}(G)$ is an \textit{optimal orientation} of $G$ if $d(D)=\bar{d}(G)$.
\indent\par Given any positive integers, $n, p_1, p_2,\ldots, p_n$, let $K_n$ denote the complete graph of order $n$ and $K(p_1,p_2,\ldots,p_n)$ denote the complete $n$-partite graph having $p_i$ vertices in the $i$th partite set for $i=1,2,\ldots, n$, where $p_1\le p_2 \le \ldots \le p_n$. The $n$ partite sets are denoted by $V_i$, $i=1,2,\ldots, n$. i.e. $|V_i|=p_i$ for $i=1,2,\ldots, n$. Furthermore, $i_j$ denotes the $j$th vertex in $V_i$ for $i=1,2,\ldots,n$, and $j=1,2,\ldots, p_i$. Thus, $K_n\cong K(p_1,p_2,\ldots,p_n)$, where $p_1=p_2=\ldots=p_n=1$.
\indent\par For general results on orientations of graphs and digraphs, we refer the reader to a survey by Koh and Tay \cite{KKM TEG 10}. Now, we introduce some results which will be found useful for our discussion later.
\indent\par Two sets $T$ and $S$ are  \textit{independent} if $T\not\subseteq S$ and $S\not\subseteq T$. If $T$ and $S$ are independent, we may say that $S$ is independent of $T$ or $T$ is independent of $S$.
\newpage
\begin{lem}\textbf{(Sperner)}\label{Sperner}
\\Let $p$ be a positive integer and let $C$ be a collection of subsets of $N_p=\{1,2,...,p\}$ such that $S$ and $T$ are independent for any two distinct sets $S$ and $T$ in $C$. Then $|C|\le {{p}\choose{\lfloor{p/2}\rfloor}}$ with equality holding if and only if all members in C have the same size, ${\lfloor{p/2}\rfloor}$ or ${\lceil{p/2}\rceil}$.
\end{lem}
\indent\par The orientation number for a general bipartite graph was determined independently by {\v S}olt\'es \cite{LS} and Gutin \cite{GG 2}.
\begin{thm} (Solt\'es \cite{LS} and Gutin \cite{GG 2})
\label{char_K(p,q)}
\\For $q\ge p \ge 2$,
\begin{equation}
 \bar{d}(K(p,q))=\left\{
  \begin{array}{@{}ll@{}}
    3, & \text{if}\ q\le{{p}\choose{\lfloor{p/2}\rfloor}}. \nonumber\\
    4, & \text{if}\ q>{{p}\choose{\lfloor{p/2}\rfloor}}, \nonumber\\
  \end{array}\right.
\end{equation}
where ${\lfloor{x}\rfloor}$ denotes the greatest integer not exceeding the real $x$.
\end{thm}

\indent\par For general $n$-partite graphs, which includes complete tripartite graphs, the following results were obtained.

\begin{thm} (Plesnik \cite{PJ 2}, Gutin \cite{GG 2}, Koh and Tan \cite{KKM TBP 2})
\\For each integer $n\ge3$, $2\le \bar{d}(K(p_1, p_2,\ldots , p_n))\le 3$.
\end{thm}

\begin{thm}(Gutin \cite{GG 2}, Koh and Tan \cite{KKM TBP 2})
\\For each integer $n\ge3$ and $p\ge 2$, $\bar{d}(K(\overbrace{p,p,\ldots,p}^{n}))=2$.
\end{thm}

\begin{thm}(Koh and Tan \cite{KKM TBP 2})
\label{thm 1.5}
\\Let $n\ge3$ and $p_1, p_2,\ldots,p_n$ be positive integers. Denote $h=\sum\limits_{k=1}^{n}{p_i}$. If 
\begin{eqnarray}
p_i >{{h-p_i}\choose{\lfloor{(h-p_i)/2\rfloor}}}\text{ for some }i=1,2, \ldots,n,\nonumber
\end{eqnarray}
then $\bar{d}(K(p_1, p_2,... , p_n))=3$.
\end{thm}

\indent\par Next, we state some existing results on complete tripartite graphs, most of which were established by Rajasekaran and Sampathkumar.
\begin{thm}(Rajasekaran and Sampathkumar \cite{RG SR})
\\For $q\ge p \ge 2$, $\bar{d}(K(1,p,q))=3$.
\end{thm}
\begin{thm}(Koh and Tan \cite{KKM TBP 3})
\label{d_bar(K(2,p,q))=2}
\\For $q\ge p \ge 2$, if $q\le{{p}\choose{\lfloor{p/2}\rfloor}}$, then $\bar{d}(K(2,p,q))=2$.
\end{thm}
\begin{thm}(Rajasekaran and Sampathkumar \cite{RG SR})
\label{d_bar(K(2,2,q))=3}
\\For $q\ge3$, $\bar{d}(K(2,2,q))=3$.
\end{thm}
\begin{thm}(Rajasekaran and Sampathkumar \cite{RG SR})
\label{d_bar(K(2,3,q))=3}
\\For $q\ge4$, $\bar{d}(K(2,3,q))=3$.
\end{thm}
\begin{thm}(Rajasekaran and Sampathkumar \cite{RG SR})
\label{thm1.10}
\\For $p \ge 4$, $4\le q \le 2p$, $\bar{d}(K(p,p,q))=2$.
\end{thm}
\indent\par Now, we proceed to further investigate the orientation number of complete tripartite graphs.
\section{A conjecture on $K(2,p,q)$}
\indent\par Based on Theorems \ref{d_bar(K(2,2,q))=3}, \ref{d_bar(K(2,3,q))=3} and an unpublished paper ``The orientation number of the complete tripartite graph $K(2,4,p)$", Rajasekaran and Sampathkumar conjectured that the converse of Theorem \ref{d_bar(K(2,p,q))=2} holds for complete tripartite graphs $K(2,p,q)$, $q\ge p \ge 5$. Ng \cite{NKL 2} showed for $q\ge p$, $\bar{d}(K(1,1,p,q))=2$ implies $q\le{{p}\choose{\lfloor{p/2}\rfloor}}$. Since an orientation $D$ of $K(2,p,q)$, where $d(D)=2$, is a spanning subdigraph of $K(1,1,p,q)$, the conjecture follows from Ng's result. In this section, we provide a different and shorter proof of the conjecture. We start with some observations which will be used in our proof later.
\begin{lem}
\label{lem1}
Let $G=K(p_1, p_2,\ldots,p_n)$, $n\ge 3$, and $D$ be an orientation of $G$. Suppose there exist vertices $i_s$ and $j_t$ for some $i, j, s$ and $t$, where $i\neq j$, $1\le i,j \le n$, $1\le s\le p_i$ and $1\le t\le p_j$, such that $O(i_s)\cap (V(G)-V_j)=O(j_t)\cap (V(G)-V_i)$. Then, $d(D)\ge 3$.
\end{lem}
\textit{Proof}: WLOG, we assume $j_t \rightarrow i_s$. It follows that $d_D(i_s, j_t)>2$ and $d(D)\ge3$.
\begin{flushright}
$\Box$
\end{flushright}

\begin{lem}
\label{lem2}
Let $D$ be an orientation of a graph $G$. Let $\tilde{D}$ be the orientation of $G$ such that $(u,v)\in E(\tilde{D})$ if and only if $(v,u) \in E(D)$. Then, $d(\tilde{D})=d(D)$.
\end{lem}
\textit{Proof}: Suppose not. Then, there exists vertices $u,v\in V(\tilde{D})$ such that $d_{\tilde{D}}(u,v)>d(D)$. Since $d_D(v,u)=d_{\tilde{D}}(u,v)$, it follows that $d_{D}(v,u)>d(D)$, yielding a contradiction.
\begin{flushright}
$\Box$
\end{flushright}

\begin{thm}
\label{char K(2,p,q)}
For any integers $q\ge p\ge 3$, if $\bar{d}(K(2,p,q))=2$, then $q\le {{p}\choose{\lfloor{p/2}\rfloor}}$.
\end{thm}
\textit{Proof}: 
\indent\par Since $\bar{d}(K(2,p,q))=2$, there exists an orientation $D$ of $K(2,p,q)$ such that $d(D)=2$.
\\
\\Case 1. $V_1 \rightarrow V_2$.
\indent\par It follows from $d_D(3_i,1_j)\le 2$, for every $i=1,2,\ldots,q$, and $j=1,2$, that $V_3 \rightarrow V_1$. Also, since $d_D(2_i, 3_j)\le 2$ for every $i=1,2,\ldots, p$, and $j=1,2,\ldots, q$, we have $V_2 \rightarrow V_3$. However, $d_D(3_i, 3_j)\ge 3$ for any $1\le i,j \le q$, $i\neq j$, which contradicts $d(D)=2$.
\indent\par Similarly, from Lemma \ref{lem2}, we cannot have $V_2\rightarrow V_1$.
\\
\\Case 2. $1_i\rightarrow V_2\rightarrow 1_{3-i}$ for exactly one of $i=1,2$.
\indent\par WLOG, we may assume that $1_1\rightarrow V_2\rightarrow 1_2$. It follows from $d_D(1_2, 3_i) \le 2$ and $d_D(3_i,1_1)\le 2$ for every $i=1,2,\ldots, q$ that $1_2\rightarrow V_3 \rightarrow 1_1$. Now, for any $i\neq j$, $1\le i, j \le q$, $d_D(3_i, 3_j)\le 2$ and thus, $O(3_i)\cap V_2$ and $O(3_j)\cap V_2$ are independent. By Sperner's Lemma, $q\le {{p}\choose{\lfloor{p/2}\rfloor}}$.
\\
\\Case 3. $1_i\rightarrow V_2$ for exactly one of $i=1,2$.
\indent\par WLOG, let $i=1$. Furthermore, we assume that $\emptyset \neq O(1_2)\cap V_2 \subset V_2$ in view of Cases 1 and 2. Hence, let $|O(1_2)\cap V_2|=k$, where $0<k<p$. Since $d_D(u,3_j)\le 2$ for every $u\in O(1_2)\cap V_2$ and every $j=1,2,\ldots, q$, it follows that $O(1_2)\cap V_2 \rightarrow V_3$. It also follows from $d_D(3_j,1_1)\le 2$ for every $j=1,2,\ldots, q$, that $V_3 \rightarrow 1_1$.
\indent\par Partition $V_3$ into $L_1$ and $L_2$ such that $L_1:=\{v\in V_3 |\text{ }v\rightarrow 1_2\}$ and $L_2:=\{v\in V_3 |\text{ }1_2\rightarrow v\}$. Note that $L_1\rightarrow V_1$. Since $d_D(2_j,v)\le 2$ for all $j=1,2,\ldots,p$, and $v\in L_1$, we have $V_2 \rightarrow L_1$. Thus, $|L_1|\le 1$, otherwise if $u,v\in L_1$, then $d_D(u,v)\ge 3$. Also, $|L_2| \le {{p-k}\choose{\lfloor{(p-k)/2}\rfloor}}$. Otherwise, by Sperner's Lemma, there exist $3_i,3_j\in L_2$ such that $O(3_i)\cap V_2 \subseteq O(3_j)\cap V_2$ for some $i\neq j$ and $1\le i, j \le q$, which implies $d_D(3_i,3_j)>2$. Hence, $q=|V_3|=|L_1|+|L_2|\le1+{{p-k}\choose{\lfloor{(p-k)/2}\rfloor}} \le1+{{p-1}\choose{\lfloor{(p-1)/2}\rfloor}}\le{{p}\choose{\lfloor{p/2}\rfloor}}$.
\indent\par Similarly, the case where $V_2\rightarrow 1_i$ for exactly one of $i=1,2$ follows from Lemma \ref{lem2}.
\\
\\Case 4. $\emptyset\neq O(1_i)\cap V_2 \subset V_2$ for $i=1,2$.
\indent\par Partition $V_3$ into the sets $L_A:=\{v\in V_3|\text{ }A\rightarrow v\rightarrow (V_1-A)\}$, where $A\subseteq V_1$. Similarly, partition $V_2$ into the sets $K_B:=\{v\in V_2|\text{ }B\rightarrow v\rightarrow (V_1-B)\}$, where $B\subseteq V_1$.
\indent\par Since $d_D(u,2_j)\le 2$ for any $u\in L_A$ and $j=1,2,\ldots,p$, it follows that $L_\emptyset\rightarrow K_\emptyset$, $L_{\{1_1\}}\rightarrow K_{\{1_1\}}\cup K_\emptyset$, $L_{\{1_2\}}\rightarrow K_{\{1_2\}}\cup K_\emptyset$ and $L_{V_1}\rightarrow V_2$. Similarly, since $d_D(u,3_j)\le 2$ for any $u\in K_B$ and $j=1,2,\ldots,q$, it follows that $K_\emptyset\rightarrow L_\emptyset$, $K_{\{1_1\}}\rightarrow L_{\{1_1\}}\cup L_\emptyset$, $K_{\{1_2\}}\rightarrow L_{\{1_2\}}\cup L_\emptyset$ and $K_{V_1}\rightarrow V_2$.
\indent\par Invoking Sperner's Lemma on each $L_A$, $A\subseteq V_1$, we have $|L_\emptyset|\le 1$, $|L_{\{1_1\}}|\le {{|K_{\{1_2\}}|}\choose{\lfloor{|K_{\{1_2\}}|/2}\rfloor}}$, $|L_{\{1_2\}}|\le {{|K_{\{1_1\}}|}\choose{\lfloor{|K_{\{1_1\}}|/2}\rfloor}}$ and $|L_{V_1}|\le 1$. Otherwise, there exist $3_i,3_j\in L_A$ such that $O(3_i)\cap V_2 \subseteq O(3_j)\cap V_2$ for some $i\neq j$ and $1\le i, j \le q$, implying $d_D(3_i,3_j)>2$. 
\\
\\Subcase 4.1. $|K_{V_1}|=0$.
\indent\par For $i=1,2$, $K_{\{1_i\}}\neq \emptyset$, since $O(1_i)\cap V_2\neq \emptyset$ by assumption. From Lemma \ref{lem1}, it follows that $L_{\{1_1\}}=L_{\{1_2\}}=\emptyset$. So, $q=|V_3|=|L_\emptyset|+|L_{\{1_1\}}|+|L_{\{1_2\}}|+|L_{V_1}|\le 1+0+0+1<{{p}\choose{\lfloor{p/2}\rfloor}}$.
\\
\\Subcase 4.2. $|K_{V_1}|>0$.
\indent\par Then, $L_{V_1}=\emptyset$ by Lemma \ref{lem1}. Recall that $|K_\emptyset|+|K_{\{1_1\}}|+|K_{\{1_2\}}|+|K_{V_1}|=p$. By Lemma \ref{lem1}, for each $i=1,2$, if $K_{\{1_i\}}\neq \emptyset$, then $L_{\{1_i\}}=\emptyset$. Hence, if $K_{\{1_1\}}\neq \emptyset$ and $K_{\{1_2\}}\neq \emptyset$, then $q=|V_3|=|L_\emptyset|+|L_{\{1_1\}}|+|L_{\{1_2\}}|+|L_{V_1}|\le 1+0+0+0$. If $K_{\{1_1\}}= \emptyset$ and $K_{\{1_2\}}\neq \emptyset$ , then $q=|L_\emptyset|+|L_{\{1_1\}}|+|L_{\{1_2\}}|\le 1+{{|K_{\{1_2\}}|}\choose{\lfloor{|K_{\{1_2\}}|/2}\rfloor}}+0 \le 1+{{p-1}\choose{\lfloor{(p-1)/2}\rfloor}}$. By symmetry, if $K_{\{1_1\}}\neq \emptyset$ and $K_{\{1_2\}}= \emptyset$, it also follows that $q\le 1+{{p-1}\choose{\lfloor{(p-1)/2}\rfloor}}$. Lastly, if $K_{\{1_1\}}=K_{\{1_2\}}=\emptyset$, it follows that that $q=|L_\emptyset|+|L_{\{1_1\}}|+|L_{\{1_2\}}|\le 1+1+1$. Therefore, $q\le{max}\big{\{}1+{{p-1}\choose{\lfloor{(p-1)/2}\rfloor}},3 \big{\}}\le{{p}\choose{\lfloor{p/2}\rfloor}}$.
\begin{flushright}
$\Box$
\end{flushright}
\begin{cor}
For any integers $p \ge 4$ and ${{p}\choose{\lfloor{p/2}\rfloor}}\ge q> 1+{{p-1}\choose{\lfloor{(p-1)/2}\rfloor}}$, let $D$ be an optimal orientation of $K(2,p,q)$, where $d(D)=2$. Then, in $D$,
\\(i) $1_i\rightarrow V_2\rightarrow 1_{3-i}\rightarrow V_3 \rightarrow 1_i$ for exactly one of $i=1,2$.
\\(ii) $\{O(3_i)\cap V_2 |\text{ }i=1,2,\ldots, q\}$ is a family of independent subsets of $\{2_1, 2_2, \ldots, 2_p\}$.
\\In particular, there are at most two optimal orientations (up to isomorphism) in the case where $q={{p}\choose{\lfloor{p/2}\rfloor}}$.
\end{cor}
\textit{Proof}: 
\indent\par Case 1 of the proof of Theorem \ref{char K(2,p,q)} shows that it is impossible for $V_1\rightarrow V_2$ or $V_2\rightarrow V_1$. Since $q> 1+{{p-1}\choose{\lfloor{(p-1)/2}\rfloor}}$ and $p\ge 4$, Cases 3 and 4 are also impossible. This leaves us with the result of Case 2, i.e. $1_i\rightarrow V_2\rightarrow 1_{3-i}\rightarrow V_3 \rightarrow 1_i$ for exactly one of $i=1,2$.
\indent\par Now, for any $i,j$ where $i\neq j$ and $1\le i,j\le q$, $3_i,3_j \in V_3$, $d(3_i, 3_j)=2$ if and only if $O(3_i)\cap V_2 \not\subseteq O(3_j)\cap V_2$. Thus, (ii) follows.
\indent\par Furthermore, if $q={{p}\choose{\lfloor{p/2}\rfloor}}$, then $|O(3_i)\cap V_2|={\lfloor\frac{p}{2}\rfloor}$ or ${\lceil\frac{p}{2}\rceil}$ by Sperner's Lemma. Thus, there are at most two optimal orientations (up to isomorphism) $D$.
\begin{flushright}
$\Box$
\end{flushright}

\indent\par Theorem \ref{char K(2,p,q)} completes the characterizaion of complete tripartite graphs $K(2,p,q)$ with $\bar{d}(K(2,p,q))=2$. Together with Theorems \ref{d_bar(K(2,p,q))=2} and \ref{d_bar(K(2,2,q))=3}, we have the following theorem. Interestingly, this characterisation has the same bounds for $q$ as the general bipartite graph $K(p,q)$. (See Theorem \ref{char_K(p,q)})
\begin{thm}
For any integers $q\ge p\ge 2$, $\bar{d}(K(2,p,q))=2$ if and only if $q\le{{p}\choose{\lfloor{p/2}\rfloor}}$.
\end{thm}


\section{Sufficient conditions for $\bar{d}(K(p,p,q))=2$}
\indent\par In this section, we provide some sufficient conditions on $p$ and $q$ so that $\bar{d}(K(p,p,q))=2$. Our result (see Theorem \ref{thm3.11}) improves significantly from the upper bound $2p$ of $q$ given in Theorem \ref{thm1.10}, especially when $p$ increases. We begin by solving a combinatorics problem, which will be of assistance later.
\begin{defn}
\label{defn 3.1}
Suppose $p\ge 4$ is an integer such that $p=kd$ for some $k,d \in \mathbb{Z^+}$, $1<k,d<p$. Denote a solution $(x_1,x_2,\ldots, x_{2d})^*$ if $(x_1,x_2,\ldots, x_{2d})$ satisfies
\begin{eqnarray} 
&&x_1+x_2+\ldots +x_{2d}=p, \label{eq: 3.1}\\
&&1\le x_i\le k-1 \text{, for }i=1,2,\ldots, 2d. \nonumber
\end{eqnarray}
\indent\par Define $\Phi^*(p,d):={{\sum\limits_{(x_1,x_2,\ldots,x_{2d})^*}{{{k}\choose{x_1}}{{k}\choose{x_2}}\ldots {{k}\choose{x_{2d}}}}}}$.
\end{defn}
\begin{defn}
Suppose $p\ge 4$ is an integer such that $p=kd$ for some $k,d \in \mathbb{Z^+}$, $1<k,d<p$. For any nonegative integers $i,j$, define $[i,j]$ to be the set of solutions $(x_1,x_2,\ldots, x_{2d})$ satisfying
\begin{eqnarray}
&&x_1+x_2+\ldots+x_{2d}=p,\nonumber\\
&&x_{s_m}=0, \text{ for }m=1,2,\ldots, i,\nonumber\\
&&x_{t_n}=k, \text{ for }n=1,2,\ldots, j,\text{ and }\nonumber\\
&&1\le x_r\le k-1 \text{ if } r\neq s_m\text{ and }r\neq t_n.\nonumber
\end{eqnarray}
Furthermore, we denote $\Phi(p,d,[i,j]):=\sum\limits_{(x_1,x_2,\ldots,x_{2d})\in [i,j]}{{{k}\choose{x_1}}{{k}\choose{x_2}}\ldots {{k}\choose{x_{2d}}}}$.
\end{defn}

\begin{rmk}The following may be verified easily.
\leavevmode
\\(a) $\Phi(p,d,[i,j])\ge 0$ for $0\le i,j \le d$.
\\(b) For each $[i,j]$ defined above, $0\le i,j \le d$.
\\(c) $\Phi(p,d,[d,d])={{2d}\choose{d}}$.
\\(d) $\Phi(p,d,[i,d])=\Phi(p,d,[d,i])=0$ for $0\le i\le d-1$.
\\(e) If $p$ is even, then $\Phi^*(p,\frac{p}{2})=2^p$.
\end{rmk}

\indent\par In the proof of our next proposition, we will make use of the following combinatorial identities which we quote without proof.
\begin{lem}
\label{lem 3.1}
For nonegative integers $x_i, n_i, n,k,r$, $n\ge 1$, $r\le k\le n$ and $x_i\le n_i$ for $i=1,2\ldots,r$ ,
\\(a) ${{n}\choose{k}}{{k}\choose{r}}={{n}\choose{r}}{{n-r}\choose{k-r}}$.
\\(b) ${{n}\choose{0}}-{{n}\choose{1}}+{{n}\choose{2}}-\ldots+(-1)^n{{n}\choose{n}}=0$.
\\(c) $\sum\limits_{x_1+x_2+\ldots x_{r}=p}{{{n_1}\choose{x_1}}{{n_2}\choose{x_2}}\ldots {{n_r}\choose{x_{r}}}}={{n_1+n_2+\ldots+n_r}\choose{p}}$. (Generalised Vandermonde's identity)
\end{lem} 

\begin{lem}
Suppose $p\ge 4$ is an integer such that $p=kd$ for some $k,d \in \mathbb{Z^+}$, $1<k,d<p$. 
Then, 
\begin{eqnarray}
\Phi (p,d,[i,j])=\sum\limits_{s=i}^{d}\sum\limits_{t=j}^{d}\Bigg[(-1)^{(s-i)+(t-j)}{{2d}\choose{s,t, 2d-(s+t)}}{{(2d-(s+t))k}\choose{(d-t)k}}{{s}\choose{i}}{{t}\choose{j}}\Bigg].\nonumber
\end{eqnarray}
\end{lem}
\textit{Proof}:
Let $\mu,\lambda$ be any two integers such that $i\le \mu \le d$ and $j\le \lambda \le d$. We proceed using a double counting method. Suppose $\alpha:={{{k}\choose{\bar{x}_1}}{{k}\choose{\bar{x}_2}}\ldots {{k}\choose{\bar{x}_{2d}}}}$, where $(\bar{x}_1,\bar{x}_2,\ldots, \bar{x}_{2d})$ is an element of $[\mu,\lambda]$. We shall show that each $\alpha$ contributes the same count to both sides of the equality.
\\
\\Case 1: $\mu=i$ and $\lambda=j$.
\indent\par On the left side, $\alpha$ is counted exactly once. The expression ${{2d}\choose{s,t,2d-(s+t)}}{{k}\choose{0}}^s{{k}\choose{k}}^t{{(2d-(s+t))k}\choose{(d-t)k}}$ represents choosing $s$ and $t$ groups from all $2d$ groups of $k$ elements to select $0$ and $k$ elements, respectively, from each group, after which $(d-t)k$ elements are selected from the remaining $(2d-(s+t))k$ elements to form a total of $p=dk$ selected elements.
\indent\par Thus, on the right, $\alpha$ is counted exactly once in the first term 
\\$(-1)^{(i-i)+(j-j)}{{2d}\choose{i,j,2d-(i+j)}}{{(2d-(i+j))k}\choose{(d-j)k}}{{i}\choose{i}}{{j}\choose{j}}={{2d}\choose{i,j,2d-(i+j)}}{{k}\choose{0}}^i{{k}\choose{k}}^j{{(2d-(i+j))k}\choose{(d-j)k}}$ 
\\and contributes a zero count in the subsequent terms 
\\${{2d}\choose{s,t,2d-(s+t)}}{{(2d-(s+t))k}\choose{(d-t)k}}={{2d}\choose{s,t,2d-(s+t)}}{{k}\choose{0}}^s{{k}\choose{k}}^t{{(2d-(s+t))k}\choose{(d-t)k}}$ if $s>i$ or $t>j$. Thus, $\alpha$ is counted once on each side.
\\
\indent\par By definition of $\alpha$, $\alpha$ is counted by the term, $\Phi (p,d,[i,j])$, on the left if and only if $[\mu,\lambda]=[i,j]$. Therefore, $\alpha$ has a zero count on the left side for the following three cases. It suffices to show that $\alpha$ contributes to a count of zero on the right in each of the following cases as well.
\\
\\Case 2: $\mu=i$ and $\lambda>j$.
\indent\par Similar to above, on the right, $\alpha$ is counted
\begin{eqnarray}
{{\lambda}\choose{j}}&\text{ times in }&{{2d}\choose{i,j,2d-(i+j)}}{{(2d-(i+j))k}\choose{(d-j)k}}\nonumber\\
{{\lambda}\choose{j+1}}&\text{ times in }&{{2d}\choose{i,j+1,2d-(i+j+1)}}{{(2d-(i+j+1))k}\choose{(d-(j+1))k}}\nonumber\\
&\vdots&\nonumber\\
{{\lambda}\choose{\lambda}}&\text{ times in }&{{2d}\choose{i,\lambda, 2d-(i+\lambda)}}{{(2d-(i+\lambda))k}\choose{(d-\lambda)k}}\nonumber
\end{eqnarray}
and none in the subsequent terms ${{2d}\choose{s,t,2d-(s+t)}}{{(2d-(s+t))k}\choose{(d-t)k}}$ if $s>i$ or $t>\lambda$.
So, $\alpha$ has a total count of $\sum\limits_{s=i}^{i}\sum\limits_{t=j}^{\lambda}[(-1)^{(s-i)+(t-j)}{{\lambda}\choose{t}}{{s}\choose{i}}{{t}\choose{j}}]=(-1)^{(i-i)}{{i}\choose{i}}\sum\limits_{t=j}^{\lambda}(-1)^{(t-j)}{{\lambda}\choose{t}}{{t}\choose{j}}$
\\$=\sum\limits_{t=j}^{\lambda}(-1)^{(t-j)}{{\lambda}\choose{j}}{{\lambda-j}\choose{t-j}}
={{\lambda}\choose{j}}\sum\limits_{t=j}^{\lambda}(-1)^{(t-j)}{{\lambda-j}\choose{t-j}}={{\lambda}\choose{j}}(0)=0$, where Lemmas \ref{lem 3.1}(a) and \ref{lem 3.1}(b) were invoked in the second and fourth equalities respectively. Thus, $\alpha$ has a zero count on each side.
\\
\\Case 3: $\mu>i$ and $\lambda=j$.
\indent\par Similar to Case 2.
\\
\\Case 4: $\mu>i$ and $\lambda>j$.
\indent\par On the right, $\alpha$ is counted ${{\mu}\choose{s}}{{\lambda}\choose{t}}$ times in the term ${{2d}\choose{s,t,2d-(s+t)}}{{(2d-(s+t))k}\choose{(d-t)k}}$, $i\le s\le \mu$ and $j\le t\le \lambda$ and 0 times if $\mu<s\le d$ or $\lambda<t\le d$. In other words, on the right, $\alpha$ is counted
\begin{eqnarray}
&&\sum\limits_{s=i}^{\mu}\sum\limits_{t=j}^{\lambda}[(-1)^{(s-i)+(t-j)}{{\mu}\choose{s}}{{\lambda}\choose{t}}{{s}\choose{i}}{{t}\choose{j}}\nonumber\\
&=&\sum\limits_{s=i}^{\mu}\Big\{(-1)^{(s-i)}{{\mu}\choose{s}}{{s}\choose{i}}\sum\limits_{t=j}^{\lambda}[(-1)^{(t-j)}{{\lambda}\choose{t}}{{t}\choose{j}}]\Big\}\nonumber\\
&=&\sum\limits_{s=i}^{\mu}\Big\{(-1)^{(s-i)}{{\mu}\choose{s}}{{s}\choose{i}}\sum\limits_{t=j}^{\lambda}[(-1)^{(t-j)}{{\lambda}\choose{j}}{{\lambda-j}\choose{t-j}}]\Big\}\nonumber\\
&=&\sum\limits_{s=i}^{\mu}\Big\{(-1)^{(s-i)}{{\mu}\choose{s}}{{s}\choose{i}}{{\lambda}\choose{j}}\sum\limits_{t=j}^{\lambda}[(-1)^{(t-j)}{{\lambda-j}\choose{t-j}}]\Big\}\nonumber\\
&=&\sum\limits_{s=i}^{\mu}\Big\{(-1)^{(s-i)}{{\mu}\choose{s}}{{s}\choose{i}}{{\lambda}\choose{j}}(0)\Big\}\nonumber\\
&=&0.\nonumber
\end{eqnarray}
times, where Lemmas \ref{lem 3.1}(a) and \ref{lem 3.1}(b) were invoked in the second and fourth equalities above respectively. Thus, $\alpha$ contributes a count of zero on each side.
\begin{flushright}
$\Box$
\end{flushright}
\begin{cor}Suppose $p\ge 4$ is an integer such that $p=kd$ for some $k,d \in \mathbb{Z^+}$, $1<k,d<p$. 
Then, 
\\(i) $\Phi^*(p,d)=\sum\limits_{s=0}^{d}\sum\limits_{t=0}^{d}[(-1)^{(s+t)}{{2d}\choose{s,t,2d-(s+t)}}{{(2d-(s+t))k}\choose{(d-t)k}}]$.
\\(ii) ${{2p}\choose{p}}=\sum\limits_{i=0}^{d}\sum\limits_{j=0}^{d}\sum\limits_{s=i}^{d}\sum\limits_{t=j}^{d}[(-1)^{(s-i)+(t-j)}{{2d}\choose{s,t,2d-(s+t)}}{{(2d-(s+t))k}\choose{(d-t)k}}{{s}\choose{i}}{{t}\choose{j}}]$.
\\(iii) $\Phi(p,d,[i,j])=\Phi(p,d,[j,i])$ for $0\le i,j \le d$.
\end{cor}
\textit{Proof}:
\\(i) This follows from the fact that  $\Phi^*(p,d)=\Phi(p,d,[0,0])$.
\\(ii) By generalised Vandermonde's identity, ${{2p}\choose{p}}=\sum\limits_{i=0}^{d}\sum\limits_{j=0}^{d}\Phi (p,d,[i,j])$.
\\(iii) Since ${{2d}\choose{s,t,2d-(s+t)}}={{2d}\choose{t,s,2d-(s+t)}}$ and ${{(2d-(s+t))k}\choose{(d-t)k}}={{(2d-(s+t))k}\choose{(d-s)k}}$, it follows that $\Phi(p,d,[j,i])=\sum\limits_{s=j}^{d}\sum\limits_{t=i}^{d}[(-1)^{(s-j)+(t-i)}{{2d}\choose{s,t,2d-(s+t)}}{{(2d-(s+t))k}\choose{(d-t)k}}{{s}\choose{j}}{{t}\choose{i}}]$
\\$=\sum\limits_{t=i}^{d}\sum\limits_{s=j}^{d}[(-1)^{(t-i)+(s-j)}{{2d}\choose{t,s,2d-(s+t)}}{{(2d-(s+t))k}\choose{(d-s)k}}{{t}\choose{i}}{{s}\choose{j}}]=\Phi(p,d,[i,j])$.
\begin{flushright}
$\Box$
\end{flushright}

\indent\par Now, we shall construct an orientation $F$ of $K(p,p,q)$, which resembles the definition of $\Phi^*(p,d)$ as its distinctive nature (see \eqref{eq: 3.1}) will aid in ensuring $d(F)=2$.
\begin{ppn}
\label{ppn 3.7}
Suppose $p\ge 4$ is an integer such that $p=kd$ for some $k,d \in \mathbb{Z^+}$, $1<k,d<p$. Then, $\bar{d}(K(p,p,q))=2$, if 
$2k+2 \le q \le \underset{d}{\max}\{\Phi^*(p,d)\}+2$,
where the maximum is taken over all positive divisors $d$ of $p$ satisfying $1<d<p$.
\end{ppn}
\textit{Proof}:  Partition $V_1\cup V_2$ into $X_1,X_2,\ldots,X_{2d}$ where 
\begin{eqnarray}
&&X_s=\{1_j| j \equiv s \text{ }(mod\text{ }d)\}, \nonumber\\
&&X_{d+s}=\{2_{(s-1)k+1},2_{(s-1)k+2},\ldots,2_{(s-1)k+k}\},\nonumber
\end{eqnarray}
for $s =1,2,\ldots,d$. Observe that $|X_r|=k$ for all $r=1,2,\ldots, 2d$. First, we define an orientation $F$ for $K(p,p,2k+2)$ as follows. (See Figure \ref{fig: 1} for $F$ when $d=3$, and $k=2$.)
\\(i) $(V_2-X_{d+s}) \rightarrow X_s\rightarrow X_{d+s}$, for $s =1,2,\ldots,d$.
\\(ii) $V_1\rightarrow 3_{2k+1}\rightarrow V_2\rightarrow 3_{2k+2}\rightarrow V_1$.
\\(iii) For $t=1,\ldots,k$,
\begin{eqnarray}
&\bullet&\{2_{k},2_{2k},\ldots,2_{dk}\}\cup (V_1-\{1_{(t-1)d+1},1_{(t-1)d+2},\ldots, 1_{(t-1)d+d}\})\rightarrow 3_t\rightarrow \nonumber\\
&& \{1_{(t-1)d+1}, 1_{(t-1)d+2},\ldots, 1_{(t-1)d+d}\}\cup (V_2-\{2_{k},2_{2k},\ldots,2_{dk}\}),\nonumber\\
&\bullet&\{1_{1},1_{2},\ldots,1_{d}\}\cup(V_2-\{2_{t},2_{t+k},\ldots, 2_{t+(d-1)k}\})\rightarrow 3_{t+k} \rightarrow \nonumber\\
&& \{2_{t},2_{t+k},\ldots, 2_{t+(d-1)k}\}\cup (V_1-\{1_{1},1_{2},\ldots,1_{d}\}).\nonumber
\end{eqnarray}
\indent\par Now, consider the case where $q>p+2$. Let $x_i=|O(3_j)\cap X_i|$ for some $j$, where $2k+2< j\le q$, and $i=1,2,\ldots, 2d$. So, for each solution $(x_1,x_2,\ldots,x_{2d})^*$ of \eqref{eq: 3.1}, there are ${{{k}\choose{x_1}}{{k}\choose{x_2}}\ldots {{k}\choose{x_{2d}}}}$ ways to choose $p$ vertices (as the outset of a vertex $3_j$), where $x_i$ vertices are selected from the set $X_i$, satisfying $1\le x_i\le k-1,$ for $i=1,2,\ldots, 2d$ and $x_1+x_2+\ldots+x_{2d}=p$. Summing over all possible solutions $(x_1,x_2,\ldots,x_{2d})^*$, there is a total of $\Phi^*(p,d):={\sum\limits_{(x_1,x_2,\ldots,x_{2d})^*}{{{k}\choose{x_1}}{{k}\choose{x_2}}\ldots {{k}\choose{x_{2d}}}}}$ of such combinations of $p$ vertices. Denote this set of combinations as $\Psi$.
\indent\par Note from (iii) that the $2k$ outsets of $3_1, 3_2,\ldots, 3_{2k}$ are elements of $\Psi$. That leaves $|\Psi|-2k=\Phi^*(p,d)-2k$ combinations of $p$ vertices of $V_1 \cup V_2$. Hence, for $2k+2<j\le q\le \underset{d}{\max} \{\Phi^*(p,d)\}+2$, we extend the definition of the above orientation so that the outset of vertices $3_{2k+3}, 3_{2k+4},\ldots, 3_q$ are these remaining elements of $\Psi$.

\begin{center}
\begin{tikzpicture}
\draw (5,9) node[circle, draw, fill=black!100, inner sep=0pt, minimum width=5pt, label=0:{$1_1$}](1_1){};
\draw (5,7.5) node[circle, draw, fill=black!100, inner sep=0pt, minimum width=5pt, label=0:{$1_4$}](1_4){};
\draw (5,6) node[circle, draw, fill=black!100, inner sep=0pt, minimum width=5pt, label=0:{$1_2$}](1_2){};
\draw (5,4.5) node[circle, draw, fill=black!100, inner sep=0pt, minimum width=5pt, label=0:{$1_5$}](1_5){};
\draw (5,3) node[circle, draw, fill=black!100, inner sep=0pt, minimum width=5pt, label=0:{$1_3$}](1_3){};
\draw (5,1.5) node[circle, draw, fill=black!100, inner sep=0pt, minimum width=5pt, label=0:{$1_6$}](1_6){};
\node[draw, inner xsep=6.5mm,inner ysep=6.5mm, fit=(1_1)(1_4), label={0:$X_1$}](x_1){};
\node[draw, inner xsep=6.5mm,inner ysep=6.5mm, fit=(1_2)(1_5), label={0:$X_2$}](x_2){};
\node[draw, inner xsep=6.5mm,inner ysep=6.5mm, fit=(1_3)(1_6), label={0:$X_3$}](x_3){};
\node[fit=(x_1)(x_2)(x_3),label={90:$V_1$}](v_1){};

\draw (-5,1.5) node[circle, draw, fill=black!100, inner sep=0pt, minimum width=5pt, label=180:{$2_6$}](2_6){};
\draw (-5,3) node[circle, draw, fill=black!100, inner sep=0pt, minimum width=5pt, label=180:{$2_5$}](2_5){};
\draw (-5,4.5) node[circle, draw, fill=black!100, inner sep=0pt, minimum width=5pt, label=180:{$2_4$}](2_4){};
\draw (-5,6) node[circle, draw, fill=black!100, inner sep=0pt, minimum width=5pt, label=180:{$2_3$}](2_3){};
\draw (-5,7.5) node[circle, draw, fill=black!100, inner sep=0pt, minimum width=5pt, label=180:{$2_2$}](2_2){};
\draw (-5,9) node[circle, draw, fill=black!100, inner sep=0pt, minimum width=5pt, label=180:{$2_1$}](2_1){};
\node[draw, inner xsep=6.5mm,inner ysep=6.5mm, fit=(2_1)(2_2), label={180:$X_4$}](x_4){};
\node[draw, inner xsep=6.5mm,inner ysep=6.5mm, fit=(2_3)(2_4), label={180:$X_5$}](x_5){};
\node[draw, inner xsep=6.5mm,inner ysep=6.5mm, fit=(2_5)(2_6), label={180:$X_6$}](x_6){};
\node[fit=(x_4)(x_5)(x_6),label={90:$V_2$}](v_2){};

\draw (0,10.5) node[circle, draw, fill=black!100, inner sep=0pt, minimum width=5pt, label=270:{$3_1$}](3_1){};
\draw (0,9) node[circle, draw, fill=black!100, inner sep=0pt, minimum width=5pt, label=270:{$3_2$}](3_2){};
\draw (0,7.5) node[circle, draw, fill=black!100, inner sep=0pt, minimum width=5pt, label=270:{$3_3$}](3_3){};
\draw (0,6) node[circle, draw, fill=black!100, inner sep=0pt, minimum width=5pt, label=270:{$3_4$}](3_4){};
\draw (0,4.5) node[circle, draw, fill=black!100, inner sep=0pt, minimum width=5pt, label=270:{$3_5$}](3_5){};
\draw (0,3) node[circle, draw, fill=black!100, inner sep=0pt, minimum width=5pt, label=270:{$3_6$}](3_6){};
\draw (0,1.5) node[circle, draw, fill=black!100, inner sep=0pt, minimum width=5pt, label=270:{$3_7$}](3_7){};
\draw (0,0) node[circle, draw, fill=black!100, inner sep=0pt, minimum width=5pt, label=270:{$3_8$}](3_8){};
\node[fit=(3_1)(3_2)(3_3)(3_4)(3_5)(3_6)(3_7)(3_8), label={90:$V_3$}](v_3){};

\draw[line,shorten <= 0.65cm, shorten >= 0.65cm] (x_1.east) to (7,8.25) to (7,12) to (-7, 12) to (-7,8.25) to (x_4.west);
\draw[line,shorten <= 0.65cm, shorten >= 0.65cm] (x_2.east) to (7.5,5.25) to (7.5,12.5) to (-7.5, 12.5) to (-7.5,5.25) to (x_5.west);
\draw[line,shorten <= 0.65cm, shorten >= 0.65cm] (x_3.east) to (7.5,2.25) to (7.5,-1) to (-7.5, -1) to (-7.5,2.25) to (x_6.west);

\draw[line,shorten <= 0.1cm, shorten >= 0.1cm] (3_5) -- (2_1);
\draw[line,shorten <= 0.1cm, shorten >= 0.1cm] (3_5) -- (2_2);
\draw[line,shorten <= 0.1cm, shorten >= 0.1cm] (3_5) -- (2_3);
\draw[line,shorten <= 0.1cm, shorten >= 0.1cm] (3_5) -- (2_4);
\draw[line,shorten <= 0.1cm, shorten >= 0.1cm] (3_5) -- (2_5);
\draw[line,shorten <= 0.1cm, shorten >= 0.1cm] (3_5) -- (2_6);

\draw[line,shorten <= 0.1cm, shorten >= 0.1cm] (3_6) -- (1_1);
\draw[line,shorten <= 0.1cm, shorten >= 0.1cm] (3_6) -- (1_2);
\draw[line,shorten <= 0.1cm, shorten >= 0.1cm] (3_6) -- (1_3);
\draw[line,shorten <= 0.1cm, shorten >= 0.1cm] (3_6) -- (1_4);
\draw[line,shorten <= 0.1cm, shorten >= 0.1cm] (3_6) -- (1_5);
\draw[line,shorten <= 0.1cm, shorten >= 0.1cm] (3_6) -- (1_6);

\draw[line,shorten <= 0.1cm, shorten >= 0.1cm] (3_1) -- (1_1);
\draw[line,shorten <= 0.1cm, shorten >= 0.1cm] (3_1) -- (1_2);
\draw[line,shorten <= 0.1cm, shorten >= 0.1cm] (3_1) -- (1_3);
\draw[line,shorten <= 0.1cm, shorten >= 0.1cm] (3_1) -- (2_1);
\draw[line,shorten <= 0.1cm, shorten >= 0.1cm] (3_1) -- (2_3);
\draw[line,shorten <= 0.1cm, shorten >= 0.1cm] (3_1) -- (2_5);

\draw[line,shorten <= 0.1cm, shorten >= 0.1cm] (3_2) -- (1_4);
\draw[line,shorten <= 0.1cm, shorten >= 0.1cm] (3_2) -- (1_5);
\draw[line,shorten <= 0.1cm, shorten >= 0.1cm] (3_2) -- (1_6);
\draw[line,shorten <= 0.1cm, shorten >= 0.1cm] (3_2) -- (2_1);
\draw[line,shorten <= 0.1cm, shorten >= 0.1cm] (3_2) -- (2_3);
\draw[line,shorten <= 0.1cm, shorten >= 0.1cm] (3_2) -- (2_5);

\draw[line,shorten <= 0.1cm, shorten >= 0.1cm] (3_3) -- (2_1);
\draw[line,shorten <= 0.1cm, shorten >= 0.1cm] (3_3) -- (2_3);
\draw[line,shorten <= 0.1cm, shorten >= 0.1cm] (3_3) -- (2_5);
\draw[line,shorten <= 0.1cm, shorten >= 0.1cm] (3_3) -- (1_4);
\draw[line,shorten <= 0.1cm, shorten >= 0.1cm] (3_3) -- (1_5);
\draw[line,shorten <= 0.1cm, shorten >= 0.1cm] (3_3) -- (1_6);

\draw[line,shorten <= 0.1cm, shorten >= 0.1cm] (3_4) -- (2_2);
\draw[line,shorten <= 0.1cm, shorten >= 0.1cm] (3_4) -- (2_4);
\draw[line,shorten <= 0.1cm, shorten >= 0.1cm] (3_4) -- (2_6);
\draw[line,shorten <= 0.1cm, shorten >= 0.1cm] (3_4) -- (1_4);
\draw[line,shorten <= 0.1cm, shorten >= 0.1cm] (3_4) -- (1_5);
\draw[line,shorten <= 0.1cm, shorten >= 0.1cm] (3_4) -- (1_6);

\draw[line,shorten <= 0.1cm, shorten >= 0.1cm] (3_7) -- (1_1);
\draw[line,shorten <= 0.1cm, shorten >= 0.1cm] (3_7) -- (1_2);
\draw[line,shorten <= 0.1cm, shorten >= 0.1cm] (3_7) -- (1_6);
\draw[line,shorten <= 0.1cm, shorten >= 0.1cm] (3_7) -- (2_1);
\draw[line,shorten <= 0.1cm, shorten >= 0.1cm] (3_7) -- (2_3);
\draw[line,shorten <= 0.1cm, shorten >= 0.1cm] (3_7) -- (2_5);

\draw[line,shorten <= 0.1cm, shorten >= 0.1cm] (3_8) -- (1_1);
\draw[line,shorten <= 0.1cm, shorten >= 0.1cm] (3_8) -- (1_2);
\draw[line,shorten <= 0.1cm, shorten >= 0.1cm] (3_8) -- (1_3);
\draw[line,shorten <= 0.1cm, shorten >= 0.1cm] (3_8) -- (2_1);
\draw[line,shorten <= 0.1cm, shorten >= 0.1cm] (3_8) -- (2_3);
\draw[line,shorten <= 0.1cm, shorten >= 0.1cm] (3_8) -- (2_6);
\end{tikzpicture}
\captionsetup{justification=centering}
{
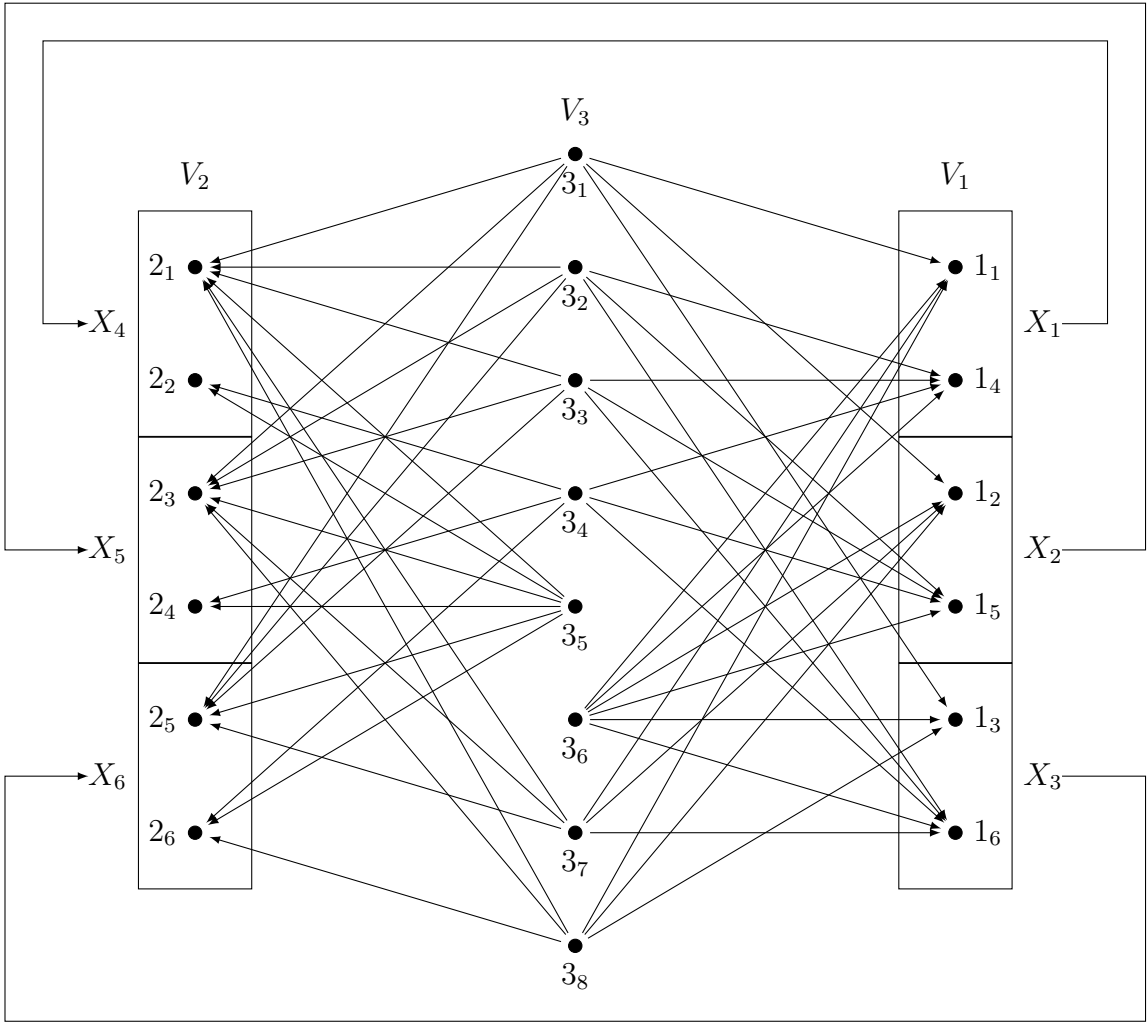
\captionof{figure}{Orientation $F$ for $d=3$, and $k=2$.\\For clarity, only the arcs from (1) $V_1$ to $V_2$ and (2) $V_3$ to $V_1$ and $V_2$ are shown.}\label{fig: 1}}
\end{center}
~\\
Claim: For all $u,v\in V(K(p,p,q))$, $d_F(u,v)\le 2$.
\\
\\Case 1: $u=1_a,v=1_b$, $a\neq b$.
\indent\par Since $1\le a,b \le p=kd$, let $a=(\alpha_1-1)d+\alpha_2$ and $b=(\beta_1-1)d+\beta_2$ for some $\alpha_i, \beta_i$, $i=1,2$, satisfying $1\le \alpha_1,\beta_1 \le k$ and $1\le \alpha_2,\beta_2 \le d$. By (iii), if $\alpha_1 \neq \beta_1$, then $1_a\rightarrow 3_{\beta_1}\rightarrow 1_b$. If $\alpha_1=\beta_1$, then it follows from $b\neq a$ that $\alpha_2\neq \beta_2$ and $a\not\equiv b\text{ }(mod$ $d)$. Therefore, $1_a\rightarrow 3_{(\alpha_2-1)k+1} \rightarrow 1_b$.
\\
\\
Case 2. $u=2_a,v=2_b$, $a\neq b$.
\indent\par Since $1\le a,b \le p=kd$, let $a=(\alpha_1-1)k+\alpha_2$ and $b=(\beta_1-1)k+\beta_2$ for some $\alpha_i, \beta_i$, $i=1,2$, satisfying $1\le \alpha_1,\beta_1 \le d$ and $1\le \alpha_2,\beta_2 \le k$. By (iii), if $\alpha_1 \neq \beta_1$, then $2_a\rightarrow 1_j \rightarrow 2_b$, where $j =d$ if $\beta_1=d$, and $\beta_1 \equiv j \text{ }(mod\text{ }d)$ otherwise. If $\alpha_1=\beta_1$, then $2_a\rightarrow 3_{\beta_2+k} \rightarrow 2_b$.
\\
\\Case 3: $u=1_a,v=2_b$.
\indent\par By (ii), $u \rightarrow  3_{2k+1} \rightarrow v$.
\\
\\
\\Case 4: $u=2_a,v=1_b$.
\indent\par By (ii), $u \rightarrow  3_{2k+2} \rightarrow v$.
\\
\\Case 5: $u=1_a,v=3_b$.
\\Subcase 5a: $b=2k+1$.
\indent\par By (ii), $V_1 \rightarrow 3_{2k+1}$.
\\
\\Subcase 5b: $b\neq 2k+1$.
\indent\par Suppose $1_a\in X_{i^*}$. Then, $1_a\rightarrow X_{d+i^*}$. Since for each $3_b$, $I(3_b)\cap X_{d+i}\neq \emptyset$ for each $i=1,2,\ldots,d$, by (ii) and (iii), let $w\in I(3_b)\cap X_{d+i^*}$. It follows that $1_a\rightarrow w \rightarrow 3_b$.
\\
\\Case 6: $u=2_a,v=3_b$.
\\Subcase 6a: $b=2k+2$.
\indent\par By (ii), $V_2 \rightarrow 3_{2k+2}$.
\\
\\Subcase 6b: $b\neq 2k+2$.
\indent\par Suppose $2_a\in X_{d+i^*}$. Then, $2_a\rightarrow X_{i}$ for all $i\neq i^*$. Since for each $3_b$, $I(3_b)\cap X_{i}\neq \emptyset$ for each $i=1,2,\ldots,d$, by (ii) and (iii), let $w\in I(3_b)\cap X_{j}$ for some $j\neq i^*$. It follows that $2_a\rightarrow w \rightarrow 3_b$.
\\
\\Case 7: $u=3_a, v=1_b$.
\\Subcase 7a: $a= 2k+2$.
\indent\par By (ii), $3_{2k+2}\rightarrow V_1$.
\\
\\Subcase 7b: $a\neq 2k+2$.
\indent\par Suppose $1_b\in X_{i^*}$. Then, $X_{d+j}\rightarrow 1_b$ for some $j\neq i^*$. Since for each $3_a$, $O(3_a)\cap X_{d+i}\neq \emptyset$ for each $i=1,2,\ldots,d$, by (ii) and (iii), let $w\in O(3_b)\cap X_{d+j}$. It follows that $3_a\rightarrow w \rightarrow 1_b$.
\\
\\Case 8: $u=3_a,v=2_b$.
\\Subcase 8a: $a= 2k+1$.
\indent\par By (ii), $3_{2k+1}\rightarrow V_2$.
\\
\\Subcase 8b: $a\neq 2k+1$.
\indent\par Suppose $2_b\in X_{d+i^*}$. Then, $X_{i^*}\rightarrow 2_b$. Since for each $3_a$, $O(3_a)\cap X_{i}\neq \emptyset$ for each $i=1,2,\ldots,d$, by (ii) and (iii), let $w\in O(3_b)\cap X_{i^*}$. It follows that $3_a\rightarrow w \rightarrow 2_b$.
\\
\\Case 9: $u=3_a,v=3_b$.
\\Subcase 9a: $a \neq 2k+1,2k+2$ and $b \neq 2k+1,2k+2$.
\indent\par Observe from (iii) that $|O(3_x)\cap(V_1\cup V_2)|=p$ for $x=a,b$. Furthermore, $O(3_a)\cap(V_1\cup V_2)\not \subseteq O(3_b)\cap(V_1\cup V_2)$ if $b\neq a$. Thus, there exists a vertex $w\in V_1\cup V_2$ such that $3_a\rightarrow w \rightarrow 3_b$.
\\
\\Subcase 9b: $a=2k+1$ and $b \neq 2k+1,2k+2$.
\indent\par $3_{2k+1}\rightarrow V_2$ and $I(3_b)\cap X_{d+i}\neq \emptyset$ for every $i=1,2,\ldots,d$, implies the existence of $w\in I(3_b)\cap V_2$. Hence, $3_a\rightarrow w \rightarrow 3_b$.
\\
\\Subcase 9c: $a=2k+2$ and $b \neq 2k+1,2k+2$.
\indent\par $3_{2k+2}\rightarrow V_1$ and $I(3_b)\cap X_{i}\neq \emptyset$ for every $i=1,2,\ldots,d$, implies the existence of $w\in I(3_b)\cap V_1$. Hence, $3_a\rightarrow w \rightarrow 3_b$.
\\
\\Subcase 9d: $a \neq 2k+1,2k+2$ and $b=2k+1$.
\indent\par $V_1\rightarrow 3_{2k+1}$ and $O(3_a)\cap X_{i}\neq \emptyset$ for every $i=1,2,\ldots,d$, implies the existence of $w\in O(3_a)\cap V_1$. Hence, $3_a\rightarrow w \rightarrow 3_b$.
\\
\\Subcase 9e: $a \neq 2k+1,2k+2$ and $b=2k+2$.
\indent\par $V_2\rightarrow 3_{2k+2}$ and $O(3_a)\cap X_{d+i}\neq \emptyset$ for every $i=1,2,\ldots,d$, implies the existence of $w\in O(3_a)\cap V_2$. Hence, $3_a\rightarrow w \rightarrow 3_b$.
\\
\\Subcase 9f: $a=2k+1$ and $b=2k+2$.
\indent\par By (ii), $3_{2k+1}\rightarrow V_2\rightarrow 3_{2k+2}$.
\\
\\Subcase 9g: $a=2k+2$ and $b=2k+1$.
\indent\par By (ii), $3_{2k+2}\rightarrow V_1\rightarrow 3_{2k+1}$.
\begin{flushright}
$\Box$
\end{flushright}
\indent\par Since $p$ may have different factorisations, the natural question to ask is which factor(s) $d$ of $p$ gives the best bound. Verification, using Maple \cite{Maple}, for all divisors $d$ of each composite integer $p\le 100$ shows that $\underset{d}{\max}\{\Phi^*(p,d)\}=\Phi^*(p,d_0)$ with $d_0$ being the smallest divisor of each $p$. Therefore, if $p$ is even, we define 
\begin{eqnarray}
\Phi_{even}(p)&:=&\Phi^*(p,2)\nonumber\\
&=&\sum\limits_{s=0}^{2}\sum\limits_{t=0}^{2}\Big[(-1)^{(s+t)}{{4}\choose{s,t,4-(s+t)}}{{(4-(s+t))\frac{p}{2}}\choose{(2-t)\frac{p}{2}}}\Big]\nonumber\\
&=&{{2p}\choose{p}}-8{{\frac{3p}{2}}\choose{p}}+12{{p}\choose{\frac{p}{2}}}-6. \nonumber
\end{eqnarray}
\indent\par Furthermore, we wish to extend Definition \ref{defn 3.1} and Proposition \ref{ppn 3.7} for prime numbers and $d=2$ seems to be the best candidate. Hence, we have the following generalisation, $\Phi_{odd}(p)$, for odd integers $p\ge 5$, which also provide a better bound than $\Phi(p,d_0)$ in cases where $p$ is odd and composite.

\begin{defn}
Suppose $p\ge 5$ is an odd integer. Denote a solution $(x_1,x_2,x_3,x_4)^{**}$ if $(x_1,x_2,x_3,x_4)$ satisfies 
\begin{eqnarray} 
\label{eq: 3.2}
&&x_1+x_2+ x_3 +x_{4}=p, \\
&&1\le x_i\le \lfloor\frac{p}{2}\rfloor \text{, for }i=1,2, \nonumber\\
&&1\le x_i\le \lfloor\frac{p}{2}\rfloor-1 \text{, for }i=3,4. \nonumber
\end{eqnarray}
Define $\Phi_{odd}(p):={\sum\limits_{(x_1,x_2,x_3,x_{4})^{**}}{{{\lfloor\frac{p}{2}\rfloor+1}\choose{x_1}}{{\lfloor\frac{p}{2}\rfloor+1}\choose{x_2}}{{\lfloor\frac{p}{2}\rfloor}\choose{x_3}} {{\lfloor\frac{p}{2}\rfloor}\choose{x_{4}}}}}$.
\end{defn}

\indent\par The following expression for $\Phi_{odd}(p)$ can be derived by exhausting all cases and is provided without proof.
\begin{lem}
If $p\ge 5$ is an odd integer, then $\Phi_{odd}(p)={{2p}\choose{p}}-4{{3x+2}\choose{x+1}}-4{{3x+1}\choose{x}}+2{{2x+2}\choose{x+1}}+8{{2x+1}\choose{x}}+2{{2x}\choose{x}}-4$, where $x=\lfloor\frac{p}{2}\rfloor$.
\end{lem}
\indent\par We shall now prove that $\Phi_{even}(p)$ and $\Phi_{odd}(p)$ are both greater than $\underset{3\le d<p}{\max}\{\Phi^*(p,d)\}$ for each $p\ge 4$.
\begin{ppn}
Suppose $p\ge 4$ is a composite integer and $d$ is a divisor of $p$, where $3\le d<p$.
\begin{equation}
\underset{3\le d<p}{\max}\{\Phi^*(p,d)\} < \left\{
  \begin{array}{@{}ll@{}}
    \Phi_{even}(p),& \text{if $p$ is even}, \nonumber\\
    \Phi_{odd}(p), & \text{if $p$ is odd}. \nonumber
  \end{array}\right.
\end{equation}
\end{ppn}

\textit{Proof}:
\\Case 1. $p$ is even.
\\
\\Claim 1. For any even integer $p\ge 14$ and any divisor $3\le d<p$ of $p$, ${{2p-\frac{p}{d}}\choose{p}}-8{{\frac{3p}{2}}\choose{p}}+12{{p}\choose{\frac{3p}{2}}}-6> 0$.
\begin{eqnarray}
&&{{2p-\frac{p}{d}}\choose{p}}-8{{\frac{3p}{2}}\choose{p}}+12{{p}\choose{\frac{3p}{2}}}-6 \nonumber\\
&\ge&{{2p-\frac{p}{d}}\choose{p}}-8{{\frac{3p}{2}}\choose{p}}\nonumber\\
&\ge&{{\frac{5p}{3}}\choose{p}}-8{{\frac{3p}{2}}\choose{p}}\nonumber\\
&>& 0.\nonumber
\end{eqnarray}
\indent\par The first inequality is due to $12{{p}\choose{\frac{3p}{2}}}\ge 6$, while the second inequality follows as $d\ge 3$ and $f(z):={{z}\choose{p}}$ is an increasing function for $z\ge p$. Since $f(z)$ is also strictly convex for $z\ge p$ and ${{\frac{5(13)}{3}}\choose{13}}-8{{\frac{3(13)}{2}}\choose{13}}>0$, the last inequality follows for all $p\ge 13$. So, Claim 1 follows.
\\
\indent\par Now, for each even integer $p\le 12$, we verified, using Maple, $\Phi^*(p,d)<\Phi_{even}(p)$ for all divisors $3\le d<p$ of $p$. (See Table \ref{tab1}.) Let $p\ge 14$ be an even integer. Note that $\sum\limits_{i=1}^{d}\sum\limits_{j=0}^{d}\Phi (p,d,[i,j])\ge {{k}\choose{0}}{{(2d-1)k}\choose{p}}={{2p-\frac{p}{d}}\choose{p}}$ as the expression ${{k}\choose{0}}{{(2d-1)k}\choose{p}}$ counts the number of ways such that none is selected from a (fixed) group of $k$ elements and $p$ elements are selected from the remaining $2d-1$ groups of $k$ elements. Also, recall that ${{2p}\choose{p}}=\sum\limits_{i=0}^{d}\sum\limits_{j=0}^{d}\Phi (p,d,[i,j])=\Phi (p,d,[0,0])+\sum\limits_{i=1}^{d}\sum\limits_{j=0}^{d}\Phi (p,d,[i,j])+\sum\limits_{j=1}^{d}\Phi (p,d,[0,j])$ by generalised Vandermonde's identity. It follows for each even integer $p\ge 14$ and each divisor $3\le d<p$ of $p$ that,
\begin{eqnarray}
&&{{2p}\choose{p}}-\Phi^* (p,d)\nonumber\\
&=&{{2p}\choose{p}}-\Phi (p,d,[0,0])\nonumber\\
&=&\sum\limits_{i=1}^{d}\sum\limits_{j=0}^{d}\Phi (p,d,[i,j])+\sum\limits_{j=1}^{d}\Phi (p,d,[0,j])\nonumber\\
&\ge& {{2p-\frac{p}{d}}\choose{p}}\nonumber\\
&>& 8{{\frac{3p}{2}}\choose{p}}-12{{p}\choose{\frac{3p}{2}}}+6\nonumber\\
&=&{{2p}\choose{p}}-\Phi_{even}(p), \nonumber
\end{eqnarray}
where the last inequality is due to Claim 1.
\\
\\Case 2. $p$ is odd and composite.
\\
\indent\par Denote $x:=\lfloor\frac{p}{2}\rfloor$.
\\
\\Claim 2. For any composite and odd integer $p\ge 17$ and any divisor $3\le d< p$ of $p$, ${{2p-\frac{p}{d}}\choose{p}}-4{{3x+2}\choose{x+1}}-4{{3x+1}\choose{x}}>0$.
\begin{eqnarray}
&&{{2p-\frac{p}{d}}\choose{p}}-4{{3x+2}\choose{x+1}}-4{{3x+1}\choose{x}} \nonumber\\
&=&{{2p-\frac{p}{d}}\choose{p}}-4{{3x+2}\choose{2x+1}}-4{{3x+1}\choose{2x+1}} \nonumber\\
&\ge&{{2p-\frac{p}{3}}\choose{p}}-8{{3x+2}\choose{2x+1}} \nonumber\\
&\ge&{{\frac{10x+5}{3}}\choose{2x+1}}-8{{3x+2}\choose{2x+1}} \nonumber\\
&>& 0\nonumber
\end{eqnarray}
\indent\par The first inequality is due to $d\ge 3$ and $f(z)$ is an increasing function for $z\ge p$. Since $f(z)$ is also strictly convex for $z\ge p$ and ${{\frac{10(8)+5}{3}}\choose{2(8)+1}}-8{{3(8)+2}\choose{2(8)+1}}>0$, the last inequality follows for all $x\ge 8$. Hence, Claim 2 follows.
\\
\indent\par For each composite and odd integer $p\le 15$, we verified, using Maple, $\Phi^*(p,d)<\Phi_{odd}(p)$ for all divisors $3\le d<p$ of $p$. (See Table \ref{tab1}.) Now, consider any composite and odd integer $p\ge 17$. As in Case 1, $\sum\limits_{i=1}^{d}\sum\limits_{j=0}^{d}\Phi (p,d,[i,j])\ge {{2p-\frac{p}{d}}\choose{p}}$. It follows for each composite and odd integer $p\ge 17$ and each divisor $3\le d<p$ of $p$ that,
\begin{eqnarray}
&&{{2p}\choose{p}}-\Phi^* (p,d)\nonumber\\
&=&{{2p}\choose{p}}-\Phi (p,d,[0,0])\nonumber\\
&=&\sum\limits_{i=1}^{d}\sum\limits_{j=0}^{d}\Phi (p,d,[i,j])+\sum\limits_{j=1}^{d}\Phi (p,d,[0,j])\nonumber\\
&\ge& {{2p-\frac{p}{d}}\choose{p}}\nonumber\\
&>& 4{{3x+2}\choose{x+1}}+4{{3x+1}\choose{x}} \nonumber\\
&\ge& 4{{3x+2}\choose{x+1}}+4{{3x+1}\choose{x}}-2{{2x+2}\choose{x+1}}-8{{2x+1}\choose{x}}-2{{2x}\choose{x}}+4 \nonumber\\
&=&{{2p}\choose{p}}-\Phi_{odd}(p), \nonumber
\end{eqnarray}
where the second last inequality follows from Claim 2.
\begin{flushright}
$\Box$
\end{flushright}

\indent\par In a way similar to Proposition \ref{ppn 3.7}, we can derive a sufficient condition for $\bar{d}(K(p,p,q))=2$ using $\Phi_{odd}(p)$ if $p$ is odd. For clarity, we summarise the results in the next theorem.
\begin{thm}
\label{thm3.11}
Suppose $p\ge 4$ is an integer. Then,
\begin{equation}
 \bar{d}(K(p,p,q))=2 \text{ if }\left\{
  \begin{array}{@{}ll@{}}
    p+2 \le q \le \Phi_{even}(p)+2, & \text{if $p$ is even}, \nonumber\\
    p+3 \le q \le \Phi_{odd}(p)+2, & \text{if $p$ is odd}. \nonumber\\
  \end{array}\right.
\end{equation}
\end{thm}
~\\
\begin{cor}
Suppose $n\ge 2$ and $p_i$ are positive integers for $i=1,2,\ldots, n$ such that $p_1+p_2+\ldots+p_r=p_{r+1}+p_{r+2}+\ldots+p_n=p\ge 4$ for some integer $r$. Let $G=K(p_1,p_2,\ldots, p_n, q)$. Then,
\begin{equation}
 \bar{d}(G)=2 \text{ if }\left\{
  \begin{array}{@{}ll@{}}
    p+2 \le q \le \Phi_{even}(p)+2, & \text{if $p$ is even}, \nonumber\\
    p+3 \le q \le \Phi_{odd}(p)+2, & \text{if $p$ is odd}, \nonumber\\
  \end{array}\right.
\end{equation}


\end{cor}
\textit{Proof}: Note that $G$ is a supergraph of $K(p,p,q)$ and $\bar{d}(K(p,p,q))=2$ by Theorem \ref{thm3.11}. So, there exists an orientation $D$ for $K(p,p,q)$, where $d(D)=2$. Partition $V(G)$ into three parts $\bigcup_{i=1}^r V_i$, $\bigcup_{i=r+1}^{n} V_i$ and $V_{n+1}$, and define an orientation $F$ for $G$ such that $D$ is a subdigraph of $F$, and edges not in $D$ are oriented arbitrarily.
\begin{flushright}
$\Box$
\end{flushright}

\indent\par For $x\ge 3$ and $p\ge2$, Koh and Tan \cite{KKM TBP 2} defined the function $f(x,p)$ to be the greatest integer such that $\bar{d}(K(\overbrace{p,p,\ldots,p}^{x},q))=2$. They posed the problem of determining $f(x,p)$. This looks like a very difficult problem. In this paper, we have made some progress for $x=2$ with results from Theorems \ref{thm1.10} and \ref{thm3.11}, where it follows that 

\begin{equation}
f(2,p) \ge \left\{
  \begin{array}{@{}ll@{}}
    \Phi_{even}(p)+2, & \text{if $p$ is even}, \nonumber\\
    \Phi_{odd}(p)+2, & \text{if $p$ is odd}. \nonumber\\
  \end{array}\right.
\end{equation}
\begin{center}
\begin{tabular}{ | P{1cm} | P{1cm}| P{10.5cm}| }
\hline
$\bm{p}$ & $\bm{d}$ &
\begin{equation}
\left\{
  \begin{array}{@{}ll@{}}
    \bm{\Phi_{even}(p)-\Phi^*(p,d)}, & \textbf{if $\bm{p}$ is even}, \nonumber\\
    \bm{\Phi_{odd}(p)-\Phi^*(p,d)}, & \textbf{if $\bm{p}$ is odd}. \nonumber
  \end{array}\right.
\end{equation}\\
\hline
4 & 2 & 16-16=0 \\ 
\hline\hline
6 & 2 & 486-486=0 \\ 
\hline
6 & 3 & 486-64=422 \\ 
\hline\hline
8 & 2 & 9,744-9,744=0 \\ 
\hline
8 & 4 & 9,744-256=9,488 \\ 
\hline\hline
9 & 3 & 39,400-14,580=24,820 \\ 
\hline\hline
10 & 2 & 163,750-163,750=0 \\ 
\hline
10 & 5 & 163,750-1,024=162,726 \\ 
\hline\hline
12 & 2 & 2,566,726-2,566,726=0 \\ 
\hline
12 & 3 & 2,566,726-1,580,096=986,630 \\ 
\hline
12 & 4 & 2,566,726-459,270=2,107,456 \\ 
\hline
12 & 6 & 2,566,726-4,096=2,562,630 \\ 
\hline\hline
14 & 2 & 39,227,538-39,227,538=0 \\ 
\hline
14 & 7 & 39,227,538-16,384=39,211,154 \\ 
\hline\hline
15 & 3 & 152,558,168-121,562,500=30,995,668 \\ 
\hline
15 & 5 & 152,558,168-14,880,348=137,677,820 \\ 
\hline\hline
16 & 2 & 595,351,056-595,351,056=0 \\ 
\hline
16 & 4 & 595,351,056-269,992,192=325,358,864 \\ 
\hline
16 & 8 & 595,351,056-65,536=595,285,520 \\ 
\hline\hline
18 & 2 & 9,038,224,134-9,038,224,134=0 \\ 
\hline
18 & 3 & 9,038,224,134-8,120,234,620=917,989,514 \\ 
\hline
18 & 6 & 9,038,224,134-491,051,484=8,547,172,650 \\ 
\hline
18 & 9 & 9,038,224,134-262,144=9,037,961,990 \\ 
\hline\hline
20 & 2 & 137,608,385,766-137,608,385,766=0 \\ 
\hline
20 & 4 & 137,608,385,766-95,227,343,750=42,381,042,016 \\ 
\hline
20 & 5 & 137,608,385,766-47,519,843,328=90,088,542,438 \\ 
\hline
20 & 10 & 137,608,385,766-1,048,576=137,607,337,190 \\
\hline
\end{tabular}
{\captionof{table}{Comparison of $\Phi^*(p,d)$ with $\Phi_{even}(p)$ and $\Phi_{odd}(p)$ for $4\le p\le 20$.}\label{tab1}}
\end{center}

\newpage
\section{Acknowledgement}

\indent\par We would like to express our gratitude to Dr Toh Pee Choon for the helpful discussions. The first author would like to thank the National Institute of Education, Nanyang Technological University of Singapore, for the generous support of the Nanyang Technological University Research Scholarship.


\end{document}